\documentclass[12pt,a4paper]{article}
\usepackage{amsmath,amstext,amssymb,amscd}

\def\Dscr{{\cal D}}
\def\Lscr{{\cal L}}
\def\Sscr{{\cal S}}

\def\Rset{{\mathbb R}}

\title{Condorcet domains of tiling type}

\author{
Vladimir I. Danilov
\thanks{Central Institute of Economics and
Mathematics of the RAS; 47, Nakhimovskii Prospect, 117418 Moscow, Russia;
email: danilov@cemi.rssi.ru}
 \and Alexander~V.~Karzanov
\thanks{Institute for System Analysis of the RAS; 9, Prospect 60 Let Oktyabrya, 117312
Moscow, Russia; email: sasha@cs.isa.ru}
 \and Gleb Koshevoy
\thanks{Central Institute of Economics and
Mathematics of the RAS; 47, Nakhimovskii Prospect, 117418 Moscow, Russia;
email: koshevoy@cemi.rssi.ru}}

\date{}

\begin{document}
\maketitle

\begin{abstract}
{ A Condorcet domain (CD) is a collection of linear orders on a set of
candidates satisfying the following property: for any choice of preferences of
voters from this collection, a simple majority rule does not yield cycles. We
propose a method of constructing ``large'' CDs by use of rhombus tiling
diagrams and explain that this method unifies several constructions of CDs
known earlier. Finally, we show that three conjectures on the maximal sizes of
those CDs are, in fact, equivalent and provide a counterexample to them.}
\medskip

 \end{abstract}


\section{Introduction} \label{sec:intr}

In the social choice theory, a {\em Condorcet domain} (further abbreviated as a
\emph{CD}) is a collection of linear orders on a finite set of candidates
(alternatives) such that if the voters choose their preferences to be linear
orders belonging to this collection, then a simple majority rule does not
yield cycles. For a survey, see, e.g.,~\cite{Monj}. A challenging problem in
the field is to construct CDs of ``large'' size. Several interesting methods
based on different ideas have been proposed in the literature.

Abello~\cite{Ab} constructed CDs by a method of completing a maximal chain in
the \emph{Bruhat lattice}. (For maximal chains in the Bruhat lattice and their
applications in combinatorics, see also~\cite{EG}.)
Chameni-Nembua~\cite{Ch-Nam} proved that covering distributive sublattices in the Bruhat
lattice are CDs. Fishburn~\cite{Fish} constructed CDs in the form of
``alternating schemes'', by using a clever combination of so-called ``never
conditions''. An alternating scheme of this sort is a representative of an
important class of CDs which we call {\em peak-pit} domains. Galambos and
Reiner~\cite{GR} developed an approach using the second Bruhat order. However,
each of those methods (which are briefly reviewed in the Appendix to this
paper) is rather indirect, and it may take some efforts to see that the objects
it generates are ``good CDs'' indeed.

In this paper we construct a class of inclusion-wise maximal, or
\emph{complete}, CDs by use of known planar graphical diagrams called
\emph{rhombus tilings}. Our construction and proofs are rather transparent and
the obtained CDs admit a good visualization. It should be noted that the
obtained class of CDs is essentially the same as each of above-mentioned
classes\footnote{The coincidence of the CD classes proposed by Abello and
by Galambos and Reiner was established in \cite{GR}.}. We
show that any peak-pit domain is a subdomain of a rhombus tiling CD (in
Theorem~4). As a consequence, we obtain that three conjectures posed,
respectively, by Fishburn, by Monjardet, and by Galambos and Reiner turn out to
be equivalent. Finally, a simple example that we construct disproves these
conjectures.

      \section{Linear orders and the Bruhat poset} \label{sec:Bruhat}

Let $X$ be a finite set whose elements are interpreted as \emph{alternatives}.
A {\em linear order} on $X$ is a complete transitive binary relation $<$ on
$X$. It ranges the elements of $X$, and we can encode a linear order $x_1
<\ldots<x_n$ on $X$ (where $n=|X|$) by the word $x_1 \ldots x _n$, regarding
$x_1$ as the least (or worst) alternative, $x_2$ as the next alternative, and
so on; then $x_n$ is the greatest (or best) alternative. The set of linear
orders on $X$ is denoted by $\mathcal L(X)$. If $Y\subset X$, we have a natural
restriction map $\mathcal L (X) \to \mathcal L (Y)$.

In what follows the ground set $X$ is identified with the set $[n]$ of integers
$1,\ldots,n$. We usually use Greek symbols, say, $\sigma$, for linear orders on
$[n]$, and write $i <_\sigma j$ rather than $i\, \sigma j$. The linear order
$1<2<\ldots<n$ is denoted by $\alpha $, and the reversed order
$n<(n-1)<\ldots<1$ by $\omega $.

Let $\Omega =\{(i,j)\ \colon i,j\in [n], i<j\}$. For a linear order $\sigma $,
a pair $(i,j)\in \Omega$ is called an {\em inversion} (w.r.t. $\alpha$) if $j
<_\sigma i$. The set of inversions for $\sigma $ is denoted by $Inv(\sigma )$.
In particular, $Inv(\alpha )=\emptyset $ and $Inv(\omega )=\Omega$.\medskip

\textbf{Definitions.} For linear orders $\sigma,\tau\in \mathcal L= \mathcal L
([n]) $, we write $\sigma \ll \tau $  if $Inv(\sigma )\subseteq Inv(\tau )$.
The relation $\ll $ on $\mathcal L $ is called the {\em weak Bruhat order}, and
the partially ordered set $(\mathcal L ,\ll )$ is called the {\em Bruhat
poset}. A linear order $\tau $ \emph{covers} a linear order
$\sigma $ if $Inv(\tau )$ equals $Inv(\sigma )$ plus exactly one inversion
(this is known to agree with the notion of covering in a poset). The
{\em Bruhat digraph} is formed by drawing a directed edge from $\sigma $ to
$\tau $ if and only if $\tau $ covers $\sigma $, and the underlying undirected
graph is called the {\em Bruhat graph}.\medskip

Clearly $\alpha $ and $\omega $ are the minimal and maximal elements of the
Bruhat poset. It is known that this poset is a lattice. Also $(\mathcal L,\ll)$
is the transitive closure of the Bruhat digraph. For $n=3$ this digraph is
drawn in Fig.~1.

\unitlength=.5mm \special{em:linewidth 0.4pt}
\linethickness{0.4pt}
\begin{picture}(90.00,84.00)(-85,0)
\put(60.00,5.00){\circle{11.00}} \put(85.00,25.00){\circle{11.00}}
\put(85.00,50.00){\circle{11.00}}
\put(60.00,70.00){\circle{11.00}}
\put(35.00,50.00){\circle{11.00}}
\put(35.00,25.00){\circle{11.00}}
\put(56.00,9.00){\vector(-3,2){17.00}}
\put(35.00,31.00){\vector(0,1){13.00}}
\put(39.00,54.00){\vector(3,2){17.00}}
\put(64.00,9.00){\vector(3,2){17.00}}
\put(85.00,31.00){\vector(0,1){13.00}}
\put(81.00,54.00){\vector(-3,2){17.00}} {\footnotesize
\put(60.00,5.00){\makebox(0,0)[cc]{123}}
\put(85.00,25.00){\makebox(0,0)[cc]{213}}
\put(85.00,50.00){\makebox(0,0)[cc]{231}}
\put(60.00,70.00){\makebox(0,0)[cc]{321}}
\put(35.00,50.00){\makebox(0,0)[cc]{312}}
\put(35.00,25.00){\makebox(0,0)[cc]{132}}}
\end{picture}

\begin{center}           Fig. 1. \end{center}

           \section{Condorcet domains} \label{sec:CD}

Let $\mathcal D \subseteq \mathcal L([n])$. We say that $\mathcal{D}$ is {\em
cyclic} if there exist three alternatives $i,j,k$ and three linear orders in
$\mathcal D $ whose restrictions to $\{i,j,k\}$ are of the form either
$\{ijk,jki,kij\}$ or $\{kji,jik,ikj\}$. An acyclic set $\mathcal{D}$ of linear
orders is called a {\em Condorcet domain} (CD). Such domains are important
since they admit aggregations (see, e.g., \cite{Monj}).

More precisely, consider a mapping $\nu :\mathcal D \to \mathbb Z_+$ (called a
$\mathcal D $-{\em opinion}), where $\nu (\sigma )$ is interpreted as the
number of voters that pick a linear order $\sigma$. Then $|\nu |=\sum
_{\sigma \in \mathcal D} \nu(\sigma)$ is the total number of voters. The
``social preference'' is defined to be the binary relation $sm(\nu )$ on $ [n]$
constructed by the majority rule: $i \ sm(\nu ) \ j \Longleftrightarrow $ the
number of voters which prefer $i$ to $j$ in their chosen linear orders is
strictly more than those having the opposite preference. When the relation
$sm(\nu )$ has no cycle for every $\mathcal D$-opinion $\nu$, the set $\mathcal
D$ is just a CD. (Indeed, it suffices to consider only $\mathcal D $-opinions
where the total number of voters is odd (cf.~\cite{Monj}). Then the relation
$sm(\nu)$ is complete, and the acyclicity of $\Dscr$ implies that $sm(\nu)$ is
a linear order on $[n]$. Conversely, if $\mathcal D $ is cyclic, then there
exists a $\Dscr$-opinion yielding a cycle in the ``social preference''.)
\smallskip

In the rest of this paper we consider only domains $\mathcal D\subset
\Lscr([n])$ containing both distinguished orders $\alpha $ and $\omega $ (this,
in fact, matches considerations in~\cite{Ab,Ch-Nam,Fish,GR,Monj}). We say that
$\mathcal D $ is {\em complete} if it is inclusion-wise maximal, i.e. adding to
$\mathcal D $ any new linear order would violate the acyclicity.

One can check that in case $n=3$ there are exactly four complete CDs. These
are:
\smallskip

a) the set of four orders 123, 132, 312 and 321. These orders are characterized
by the property that the alternative 2 is never the worst. We call this CD the
{\em peak domain} (for $n=3$) and denote it as $\mathcal D_3 (\cap )$.
\smallskip

b) the set of orders 123, 213, 231, 321. In these orders the alternative 2 is
never the best. This CD is called the {\em pit domain} and denoted by $\mathcal
D_3 (\cup )$. \smallskip

c) the set $\{123, 213, 312, 321\}$. Here the alternative 3 is never the
middle. We denote this domain by $\mathcal D_3 (\to)$. \smallskip

d) the set $\{123, 132, 231, 321\}$, denoted by $\mathcal D_3 (\gets)$. Here
the alternative 1 is never the middle.
\smallskip

A {\em casting} is meant to be a mapping $c$ of the set ${[n] \choose  3}$ of
triples $ijk$ ($i<j<k$) to $\{\cap ,\cup ,\to, \gets\}$. For a casting $c$, we
define $\mathcal D(c)$ to be the set of linear orders $\sigma \in \mathcal
L([n])$ whose restrictions to each triple $ijk$ (further denoted as $\sigma
|_{ijk}$) belongs to $\mathcal D_3 (c(ijk))$. The following assertions are
immediate.\medskip

{\bf Proposition 1.} (i) {\em For any casting $c$, ~$\mathcal D (c)$ is a
Condorcet domain. }

(ii) {\em Any Condorcet domain is contained in a set $\mathcal D (c)$, where
$c$ is a casting.}\medskip

Note that a casually chosen casting may produce a small and/or non-complete CD.
As Fishburn writes in \cite{Fish}: ``.. it is far from obvious how the
restrictions should be selected jointly to produce a large acyclic set.'' In
the next section we describe and examine a simple geometric construction
generating a representable class of complete CDs.

    \section{Rhombus tilings and related CDs} \label{sec:tiling-CD}

The complete CDs that we are going to introduce one-to-one correspond to
certain  geometric arrangements on the plane, called rhombus tilings. In this
section we recall this notion, review basic properties of tilings needed to us,
and finally we establish some facts about related CDs.\medskip

A. In the upper half-plane $\mathbb R\times \mathbb R_{>0}$, we fix $n$ vectors
$\xi_1 ,\ldots,\xi_ n$ going in this order clockwise around $(0,0)$ and having
the same length. The sum of segments $[0,\xi_ i]$, $i=1,\ldots,n$, forms a
\emph{zonogon}, denoted by $Z=Z_n$. This is the center-symmetric $2n$-gon
formed by the points $\sum _i a_i \xi_ i$ over all $0\le a _i\le 1$. Two
vertices of the zonogon are distinguished: the \emph{bottom} vertex
$b(Z):=(0,0)$ and the \emph{top} vertex $t(Z):=\xi_1+\ldots+\xi_n$. A rhombus
congruent to the sum of two segments $[0,\xi_i ]$ and $[0,\xi_ j]$, where $1\le
i<j\le n$, is called an $ij$-{\em tile}, or simply a {\em tile}.

A {\em rhombus tiling} (or simply a {\em tiling}) is a subdivision $T$ of the
zonogon into a set of tiles satisfying the following condition: if two tiles
intersect, then their intersection consists of a single vertex or a single
(closed) edge. The set of tiles of $T$ is denoted by $Rho(T)$. Figures 2 and 4
illustrate examples of rhombus tilings.

We associate to a tiling $T$ the planar directed graph $G_T=(V_T,E_T)$ whose
vertices and edges are those occurring in the tiles and the edges are oriented
upward. The tiles of $T$ are just the (inner two-dimensional closed) faces of
$G_T$. An edge congruent to $\xi_i$ is called an $i$-\emph{edge}, or an edge of
\emph{color} $i$.

We will need two more definitions. First, since all edges of $G_T$ are directed
upward, this digraph is acyclic and any maximal directed path in it goes from
$b(Z)$ to $t(Z)$. We call such a path a {\em snake} of $T$. In particular, the
zonogon is bounded by two snakes, namely, those forming the left boundary
$lbd(Z)$ and the right boundary $rbd(Z)$ of $Z$; note that the sequence of edge
colors in the former (latter) gives the linear order $\alpha$ (resp. $\omega$).

Second, for $i\in[n]$, we apply the term an $i$-{\em track} (borrowed
from~\cite{KS}) to a maximal alternating sequence
$Q=(e_0,F_1,e_1,\ldots,F_k,e_k)$ formed by $i$-edges and different tiles, where
$e_{j-1},e_j$ are opposite edges of a tile $F_j$ (other known names for $Q$ are
``de Bruijn line''~\cite{dB}, ``dual $i$-path'', ``$i$-stripe''.) Note that the
projections of $e_0,\ldots,e_k$ to a line orthogonal to $\xi_i$ give a monotone
sequence of points (since consecutive tiles in $Q$ do not overlap). This
implies that $Q$ is not cyclic, is determined uniquely up to reversing,
contains all $i$-edges of $T$, and connects the pair of $i$-edges on the
boundary of the zonogon. We assume for definiteness that the $i$-track begins
(with the edge $e_0$) on the left boundary of $Z_n$, and ends (with $e_k$) on
the right boundary. \medskip

B. Next we exhibit some properties of tilings. One important use of tracks
consists in the following. When removing the $i$-track $Q$ from the zonogon
(i.e. removing the interiors of the edges and tiles of $Q$), we obtain two
connected regions $L_i,U_i$ such that: $L_i$ (the \emph{lower} region) contains
the bottom vertex $b(Z)$ and $U_i$ (the \emph{upper} region) contains the top
vertex $t(Z)$; the edges of $G_T$ connecting these regions are exactly the
$i$-edges $e_0,\ldots,e_k$ and these are directed from $L_i$ to $U_i$; gluing
$L_i$ with $U_i$ shifted by $-\xi_i$ produces the $(n-1)$-zonogon $Z'$
generated by the vectors $\xi_1,\ldots,\xi_{i-1},\xi_{i+1},\ldots,\xi_n$.
Moreover, removing from $T$ the tiles of $Q$ (and shifting those in $U_i$ by
$-\xi_i$) gives a rhombus tiling $T'$ of $Z'$; we call $T'$ the
\emph{reduction} of $T$ by the color $i$ and denote it as $T|_{[n]-i}$.

Using this operation and some other simple constructions and observations, one
can demonstrate a number of rather elementary properties of tilings. Among
these, the following nice properties of $T$ are known. \medskip

{\bf Proposition~2.} {\rm(i)} {\em Any snake $S$ intersects an $i$-track at
exactly one $i$-edge. Therefore, $S$ contains exactly $n$ edges and the
sequence of edge colors along $S$ gives a linear order on $[n]$.}

{\rm(ii)} {\em For any $1\le i<j\le n$, there is exactly one $ij$-tile in $T$.
This yields a natural bijection $\psi: Rho(T)\to\Omega$ (which maps an
$ij$-tile to the pair $(i,j)\in\Omega$).}

{\rm(iii)} {\em For a snake $S$ of $T$, let $\sigma$ be the linear order
determined by $S$, and let $L(S)$, or $L(\sigma)$, denote the set of tiles of
$T$ lying on the left from $S$, i.e. those contained in the region bounded by
$S$ and $lbd(Z)$. Then $\psi(L(\sigma))=Inv(\sigma)$.}

{\rm(iv)} {\em For a snake $S$, there exist two consecutive edges $e,e'$ in $S$
(where $e$ precedes $e'$) which have colors $i$ and $j$, respectively, and
belong to a tile $\rho\in Rho(T)$ so that: {\rm (a)} if $S\ne lbd(Z)$ then
$i>j$ and $\rho$ lies on the left from $S$, and {\rm (b)} if $S\ne rbd(Z)$ then
$i<j$ and $\rho$ lies on the right from $S$.}
\medskip

\textbf{Remark.} ~These facts (or somewhat close to them) were established in
several works, possibly being formulated in different terms. See,
e.g.,~\cite{El,Good97,GR,KS,Zig}. Some authors (e.g., in~\cite{GR}) prefer to
operate in terms of so-called \emph{commutation classes} of \emph{pseudo-line
arrangements} (visualizing \emph{reduced words for permutations},
cf.~\cite{BFZ}). Such objects, related to rhombus tilings via planar duality,
are in fact equivalent to \emph{simple wiring diagrams} (a special case of
wirings studied in~\cite{DKK-10}). The latter diagram can be introduced as a
set of curves (``wires'') $\zeta_1,\ldots,\zeta_n$ in the strip $[0,1]\times
\Rset$ with the following properties: $\zeta_i$ begins at the point $(0,i)$ and
ends at the point $(1,n-i)$; any two wires intersect at exactly one point; and
no three wires have a common point. This is bijective (up to an isotopy) to a
rhombus tiling $T$ in which an $ij$-tile corresponds to the intersection point
of wires $\zeta_i,\zeta_j$ and an $i$-track corresponds to the wire $\zeta_i$.
In their turn, the snakes of $T$ correspond to the so-called \emph{cutpaths} in
the wiring (in terminology of~\cite{GR}).\medskip

In light of (i) in Proposition~2, we will not distinguish between snakes $S$
and their corresponding linear orders $\sigma$, denoting the snake as $\mathcal
S (\sigma )$ and saying that the linear order $\sigma $ is {\em compatible}
with the tiling $T$. The set of linear orders compatible with $T$ is denoted by
$\Sigma (T)$.\medskip

\textbf{Example 1.} When $n=3$, there are exactly two tilings of the zonogon
(hexagon) $Z _3$, as depicted below. Here the set $\Sigma (T)$ consists of four
orders, namely: 123, 132, 312, 321. This is precisely the peak domain $\mathcal
D (\cap )$. In its turn, the set $\Sigma (T')$ consists of four orders 123,
213, 231, 321, which is just the pit domain $\mathcal D (\cup )$.

\unitlength=.8mm \special{em:linewidth 0.4pt}
\linethickness{0.4pt}
\begin{picture}(115.00,47.00)(-10,0)
\put(40.00,5.00){\vector(-3,2){15.00}}
\put(25.00,15.00){\vector(0,1){15.00}}
\put(25.00,30.00){\vector(3,2){15.00}}
\put(40.00,6.00){\vector(-3,2){15.00}}
\put(25.00,16.00){\vector(4,3){14.00}}
\put(39.50,26.67){\vector(0,1){12.33}}
\put(40.00,5.00){\vector(3,2){15.00}}
\put(55.00,15.00){\vector(0,1){15.00}}
\put(55.00,30.00){\vector(-3,2){15.00}}
\put(40.00,6.00){\vector(3,2){15.00}}
\put(55.00,16.00){\vector(-4,3){14.00}}
\put(40.50,26.67){\vector(0,1){12.33}}

\put(100.00,5.00){\vector(-3,2){15.00}}
\put(85.00,15.00){\vector(0,1){15.00}}
\put(85.00,30.00){\vector(3,2){15.00}}
\put(100.00,5.00){\vector(3,2){15.00}}
\put(115.00,15.00){\vector(0,1){15.00}}
\put(115.00,30.00){\vector(-3,2){15.00}}
\put(100.00,40.00){\vector(0,0){0.00}}
\put(99.50,6.00){\vector(0,1){15.00}}
\put(99.00,21.00){\vector(-3,2){13.00}}
\put(86.00,29.67){\vector(3,2){14.00}}
\put(100.500,6.00){\vector(0,1){15.00}}
\put(101.00,21.00){\vector(3,2){13.00}}
\put(114.00,29.67){\vector(-3,2){14.00}}

\put(40.00,.00){\makebox(0,0)[cc]{$T$}}
\put(100.00,.00){\makebox(0,0)[cc]{$T'$}}
\end{picture}
 \begin{center}
     Fig. 2.
  \end{center}

\noindent So the domains $\Sigma (T)$ and $\Sigma (T')$ in this example are
CDs. We will explain later that a similar property holds for any rhombus
tiling. \medskip

Next, the snakes of a tiling $T$ of the zonogon $Z=Z_n$ are partially ordered
``from left to right'' in a natural way. The minimal element is the leftmost
snake $\mathcal S (\alpha )=lbd(Z)$, and the maximal element is the rightmost
snake $\mathcal S (\omega )=rbd(Z)$. The corresponding poset is a (distributive)
lattice in which for two snakes $S$ and $S'$, their greatest lower bound
$S\wedge S'$ coincides with their ``left envelope'', and the least upper bound
$S\vee S'$ coincides with the ``right envelope''. In terms of left regions of
snakes (cf.~Proposition~2(iii)), we have $L(S\wedge S')=L(S)\cap L(S')$ and
$L(S\vee S')=L(S)\cup L(S')$.

Thus, we obtain a natural partial order $\prec$ on the set $\Sigma(T)$ of
linear orders, defined by $\sigma\prec\tau \Leftrightarrow L(\sigma)\subset
L(\tau)$. Moreover, by~(iii) in Proposition~2, the relation $L(\sigma)\subset
L(\tau)$ is equivalent to $Inv(\sigma)\subset Inv(\tau)$, and therefore the
partial order $\prec$ on $\Sigma(T)$ is induced by the the weak Bruhat order
$\ll$ on $\mathcal L([n])$.

In its turn,~(iv) in Proposition~2 shows that if a snake $\mathcal S (\tau )$
lies on the right from a snake $\mathcal S (\sigma )$ and there is no snake
between them, then these snakes differ by a single tile. This leads to a
sharper version of the above property, namely: {\em the covering relations on
the poset $\Sigma(T)$ (w.r.t. $\prec$) are induced by covering relations on the
Bruhat poset}. As a consequence, we obtain the following \medskip

\textbf{Corollary 1.} {\em Any maximal chain in the poset $\Sigma (T)$ is a
maximal chain in the Bruhat poset $(\mathcal L ,\ll )$}. \medskip

C. In the rest of this section we show that for any rhombus tiling $T$ of
$Z_n$, the set $\Sigma (T)$ is a CD.

We use the track reducing operation defined above. Take the reduction
$T'=T|_{[n]-i}$ of $T$ by an alternative $i$. Then any snake $\mathcal S
(\sigma )$ compatible with $T$ is transformed into a snake corresponding to the
restricted linear order $\sigma |_{[n]-i}$ and compatible with $T'$. This gives
the restriction map
                                                 $$
                              \Sigma (T)\to \Sigma      (T|_{[n]-i}).
                                           $$
Making a sequence of reducing operations, we can reach any subset $X\subset
[n]$ and obtain the corresponding restriction map
      $$
                                \Sigma(T) \to \Sigma(T|_X).
      $$

{\bf Theorem 1.} {\em  The set $\Sigma(T)$ is a complete Condorcet
domain.}\medskip

\emph{Proof.} ~Consider the restrictions of linear orders from $\Sigma (T)$ to
a triple $ijk$. By reasonings above, they belong to $\Sigma (T|_{ijk})$. The
obtained domain is either $\mathcal D (\cup )$ or $\mathcal D (\cap )$ (defined
in Section~\ref{sec:CD}). Therefore, $\Sigma(T)$ is a CD
(cf.~Proposition~1(i)).

To check the completeness of $\Sigma(T)$, let us try to add to it a new linear
order $\rho $. Then the corresponding path $\mathcal S (\rho )$ drawn in $Z_n$
is not contained in $G_T$. Let $e$ be the \emph{first} edge of $\mathcal S
(\rho )$ which is not an edge of $T$, and let $v$ be the beginning vertex of
$e$. Then the part $P$ of $\Sscr(\rho)$ from $b(Z_n)$ to $v$ lies in $G_T$.
Three cases are possible, as depicted in Figure~3.

\unitlength=.8mm \special{em:linewidth 0.4pt}
\linethickness{0.4pt}
\begin{picture}(130.00,35)(-20,5)
\put(70.00,5.00){\vector(2,3){10.00}}
\put(70.00,5.00){\vector(0,1){15.00}}
\put(70.00,5.00){\vector(-1,1){15.00}} \put(55,20){\line(2,3){10}}
\put(65,35){\line(1,-1){15}}

\put(30,20){\circle*{2}} \put(70,5){\circle*{2}}
\put(110,20){\circle*{2}}

\put(58.00,12.00){\makebox(0,0)[cc]{$i$}}
\put(80.00,12.00){\makebox(0,0)[cc]{$k$}}
\put(73.00,15.00){\makebox(0,0)[cc]{$e$}}
\put(66,5){\makebox(0,0)[cc]{$v$}}
\put(20.00,5.00){\vector(2,3){10.00}}
\put(30.00,20.00){\vector(-2,3){10.00}}
\put(20.00,5.00){\vector(-2,3){10.00}}
\put(10.00,20.00){\vector(2,3){10.00}}
\put(30.00,20.00){\vector(-1,1){15.00}}
\put(20.00,25.00){\makebox(0,0)[cc]{$e$}}
\put(33.5,20){\makebox(0,0)[cc]{$v$}}
\put(120.00,5.00){\vector(-2,3){10.00}}
\put(110.00,20.00){\vector(2,3){10.00}}
\put(120.00,5.00){\vector(2,3){10.00}}
\put(130.00,20.00){\vector(-2,3){10.00}}
\put(110.00,20.00){\vector(4,3){20.00}}
\put(120.00,24.00){\makebox(0,0)[cc]{$e$}}
\put(107,20){\makebox(0,0)[cc]{$v$}}
\put(70.00,20.00){\vector(1,3){5.00}}

\put(11.00,12.00){\makebox(0,0)[cc]{$i$}}
\put(11.00,28.00){\makebox(0,0)[cc]{$k$}}

\put(130,28.00){\makebox(0,0)[cc]{$i$}}
\put(130,12.00){\makebox(0,0)[cc]{$k$}}
 \end{picture}
 \begin{center}           Fig. 3 \end{center}

Consider the middle case. Let the edge  $e$ have color $j$, and let the tile of
$T$ whose interior meets $e$ be an $ik$-tile $Q$. Then $i<j<k$. Clearly the
part $P$ of $\Sscr(\rho)$ cannot contain an edge with color in $\{i,j,k\}$.
Hence, in the linear order $\rho$ the alternative $j$ occurs earlier than each
of $i,k$. Two subcases are possible: either $j<_\rho i<_\rho k$ or $j<_\rho
k<_\rho i$. In the first subcase, compare $\rho$ with two linear orders from
the domain $\Sigma (T)$: a linear order $\sigma$ that follows the path $P$ and
then the left side of $Q$, yielding the relation $i<_\sigma k<_\sigma j$, and
the linear order $\omega$, yielding $k<_{\omega} j <_{\omega} i$. This gives a
cyclic triple. In the second subcase, act symmetrically, by comparing $\rho$
with a linear order $\tau$ that follows $P$ and the right side of $Q$ (yielding
$k<_\tau i<_\tau j$) and the linear order $\alpha$ (yielding $i<_\alpha
j<_\alpha k$), again obtaining a cyclic triple.

Two other cases are examined in a similar way. $\hfill\Box$ \medskip

We refer to a domain of the form  $\Sigma (T)$ as a {\em Condorcet
domain of tiling type}, or a {\em tiling CD}.

\section{Tiling CDs and peak-pit domains}

A set $\mathcal D\subset \mathcal L([n]) $ is called a {\em peak-pit domain} if
for each triple $i<j<k$ in $[n]$, the peak condition or the pit one is
satisfied (in the sense that the projection of $\Dscr$ to $\{i,j,k\}$
is contained either in the peak
domain $\mathcal D_3 (\cap )$ or in the pit domain $\mathcal D_3 (\cup )$ (with
$ijk$ in place of 123) or in both). We have the following property (cf. the proof of Theorem~1):
   \begin{itemize}
\item[($\ast$)]\qquad
\emph{any tiling CD is a peak-pit domain.}
   \end{itemize}
The converse property is valid as well.
\medskip

{\bf Theorem 2}. {\em Any peak-pit domain is contained in a tiling CD.}
\medskip

To prove this assertion (which is less trivial) we need some definitions and
preliminary observations.

Let $\sigma\in\Lscr([n])$. A subset $X\subseteq [n]$ is called an {\em ideal}
of $\sigma $ if $x\in X$ and $y<_\sigma x$ imply $y\in X$. In other words, if
$\sigma$ is represented as a word $i_1 \ldots i _n$, then an ideal of $\sigma$
corresponds to an initial segment of this word. Let $Id(\sigma )$ denote the
set of ideals of $\sigma $ (including the empty set). In particular,
$Id(\alpha)$ consists of the intervals $[0]$, $[1],\ldots,[n-1]$, $[n]$.

We associate to a collection $\mathcal D \subseteq \Lscr([n])$ the following
set-system
  $$
Id(\mathcal D )=\cup _{\sigma \in \mathcal D }Id(\sigma ).
  $$

{\bf Example 2.} ~Let  $\mathcal D $ be the peak domain for $n=3$; it consists
of four orders 123, 132, 312, and 321. Then $Id(\mathcal D )$ consists of seven
sets $\emptyset $, $\{1\}$, $\{3\}$, $\{1,2\}$, $\{1,3\}$, $\{2,3\}$, and $\{1,2,3\}$=[3], that is, of all subsets of
$[3]$ except for $\{2\}$. In its turn, for the pit domain $\mathcal D' $,
~$Id(\mathcal D' )$ consists of all subsets of $[3]$ except for $\{1,3\}$.\medskip

Consider a tiling $T$. We associate to each vertex $v$ in it a certain subset
$sp(v)$ of $[n]$, as follows. Let $\mathcal S (\sigma )$ be a snake passing
$v$. Then $sp(v)$ is the ideal of $\sigma $ corresponding to the part of
$\mathcal S (\sigma )$ from the beginning to $v$ (the set $sp(v)$ does not
depend on the choice of a snake $\sigma $ passing $v$). This is equivalent to
saying that $sp(v)$ consists of the elements $i\in [n]$ such that the $i$-track
goes below the vertex $v$ (in view of Proposition~2(i)). The collection of sets
$sp(v)$ for all vertices $v$ of $T$ is denoted by $Sp(T)$ and called the
\emph{spectrum} of $T$ (following terminology in~\cite{DKK-10}). One can check
that a linear order $\sigma $ belongs to $ \Sigma (T)$ if and only if the
inclusion $Id(\sigma )\subset Sp(T)$ holds.\medskip

{\em Proof of Theorem 2}. ~Let  $\mathcal D\subset\Lscr([n])$ be a peak-pit
domain. Our aim is to show the existence of a tiling $T$ such that $Id(\mathcal
D)\subseteq Sp(T)$. We use a criterion due to Leclerc and Zelevinsky~\cite{LZ}
on a system of subsets of $[n]$ that can be extended to the spectrum $Sp(T)$ of
a tiling $T$. (Strictly speaking, the criterion in~\cite{LZ} concerns
set-systems associated with pseudo-line arrangements (which correspond, in a
sense, to rhombus tilings, cf.~\cite{El}). For a direct proof, in terms of
tilings, see~\cite[Sec.~5.3]{DKK-08}.)

Two subsets $A,B$ of $[n]$ are said to be \emph{separated} (more precisely,
\emph{strongly separated}, in terminology of~\cite{LZ}) from each other if the
convex hulls of $A\setminus B$ and $B\setminus A$ (viz. the minimal intervals
containing these sets) are disjoint. For example, the sets $\{1,2\}$ and
$\{2,4\}$ are separated, whereas $\{1,3\}$ and $\{2\}$ are not. In particular,
$A$ and $B$ are separated if one includes the other. A collection of sets is
called {\em separated} if any two sets in it are separated.\medskip

{\bf Theorem 3 \cite{LZ}.} {\em The spectrum $Sp(T)$ of any rhombus tiling $T$
is separated. Conversely, if $\mathcal X $ is a separated set-system on $[n]$, then
there exists a tiling $T$ of $Z_n$ such that $\mathcal X \subset
Sp(T)$.}\medskip

Due to this theorem, it suffices to show that for a peak-pit domain $\mathcal D
$, the system  $Id(\mathcal D )$ is  separated. Suppose this is not so for some
$\mathcal D $. Then there are two sets $A,B\in Id(\mathcal D)$ and a triple
$i<j<k$ in $[n]$ such that $A$ contains $j$ but none of $i,k$, whereas $B$
contains $i,k$ but not $j$. Restrict the members of $\mathcal D$ to the set
$\{i,j,k\}$. Then $Id(\mathcal D |_{ijk})$ contains both sets $\{j\}$ and
$\{i,k\}$. Thus, we are neither in the peak nor in the pit domain case, as we
have seen in Example 2. \hfill$\Box$\medskip

Now we combine Theorem~2 and a slight modification of property~($\ast$) (in the
beginning of this section), yielding the main assertion in this paper. Let us
say that a domain $\mathcal{D}$ is {\em semi-connected} if the linear orders
$\alpha$ and $\omega$ can be connected in the Bruhat graph by a path in which
all vertices belong to $\mathcal{D}$.\medskip

{\bf Theorem 4.} (i) {\em Every domain of tiling type is semi-connected}.

(ii) {\em Every semi-connected Condorcet domain is a peak-pit domain.}

(iii) {\em Every peak-pit domain is contained in a domain of tiling
type.}\medskip

{\em Proof.} ~Any domain of the form $\Sigma (T)$ is semi-connected since it
contains a maximal chain of the Bruhat poset (cf. Corollary~1), yielding~(i).

It is easy to see that the semi-connectedness preserves under reducing
alternatives. Because of this, we can restrict ourselves to the case $n=3$. In
this case there exist exactly four CDs. Two of them, where one of the
alternatives 1 and 3 is never the middle, are not semi-connected. The other two
domains are semi-connected; they are just the peak and pit domains. This
implies~(ii).

Claim~(iii) is just Theorem 2. \hfill$\Box$ \medskip

As a consequence, we obtain that the CDs constructed by Abello\cite{Ab},
Chameni-Nembua~\cite{Ch-Nam}, and Galambos and Reiner~\cite{GR} (see the
Appendix for a brief outline), as well as the maximal peak-pit domains, are CDs
of tiling type. Moreover, all these classes of CDs are equal.

\section{On Fishburn's conjecture}

Fishburn \cite{Fish} constructed Condorcet domains by the following method. For
a set of linear orders and a triple $i<j<k$, Fishburn's ``never condition''
$jN1$ means the requirement that, in the restriction of each of these linear
orders to $\{i,j,k\}$, the alternative $j$ is never the worst. This is exactly
the above-mentioned ``peak condition'' for $ijk$. Similarly, the ``never
condition'' $jN3$ (saying that ``the alternative $j$ is never the best'')
coincides with the ``pit condition'' for $ijk$.

Fishburn's \emph{alternating scheme} is defined by imposing, for each triple
$i<j<k$, the peak condition when $j$ is even, and the pit condition when $j$ is
odd. The set of linear orders (individually) obeying these conditions
is called \emph{Fishburn's domain}, and its cardinality is denoted by
$\Phi (n)$.

By Theorem~2, Fishburn's domain $\mathcal D$ is contained in a CD of tiling
type. Also it is a complete CD, as is shown in~\cite{GR}. So $\mathcal D$ is a
tiling CD. The corresponding tiling for $n=8$ is drawn in Fig.~4.

\unitlength=.42mm \special{em:linewidth 0.4pt}
\linethickness{0.4pt}
\begin{picture}(120.00,95)(-110,0)
\put(60.00,5.00){\vector(-3,1){30.00}}
\put(30.00,15.00){\vector(-3,2){15.00}}
\put(15.00,25.00){\vector(-1,1){10.00}}
\put(5.00,35.00){\vector(-1,2){5.00}}
\put(0.00,45.00){\vector(1,2){5.00}}
\put(5.00,55.00){\vector(1,1){10.00}}
\put(15.00,65.00){\vector(3,2){15.00}}
\put(30.00,75.00){\vector(3,1){30.00}}
\put(60.00,5.00){\vector(3,1){30.00}}
\put(90.00,15.00){\vector(3,2){15.00}}
\put(105.00,25.00){\vector(1,1){10.00}}
\put(115.00,35.00){\vector(1,2){5.00}}
\put(120.00,45.00){\vector(-1,2){5.00}}
\put(115.00,55.00){\vector(-1,1){10.00}}
\put(105.00,65.00){\vector(-3,2){15.00}}
\put(90.00,75.00){\vector(-3,1){30.00}}
\put(60.00,5.00){\vector(-1,1){10.00}}
\put(60.00,5.00){\vector(1,2){5.00}}
\put(60.00,5.00){\vector(3,2){15.00}}
\put(75.00,15.00){\vector(3,1){30.00}}
\put(75.00,15.00){\vector(1,2){5.00}}
\put(65.00,15.00){\vector(3,2){15.00}}
\put(65.00,15.00){\vector(-1,1){10.00}}
\put(50.00,15.00){\vector(1,2){5.00}}
\put(50.00,15.00){\vector(-3,1){30.00}}
\put(30.00,15.00){\vector(-1,1){10.00}}
\put(20.00,25.00){\vector(-3,2){15.00}}
\put(20.00,25.00){\vector(1,2){5.00}}
\put(25.00,35.00){\vector(-3,2){15.00}}
\put(55.00,25.00){\vector(-3,1){30.00}}
\put(55.00,25.00){\vector(3,2){15.00}}
\put(80.00,25.00){\vector(-1,1){10.00}}
\put(80.00,25.00){\vector(3,1){30.00}}
\put(110.00,35.00){\vector(1,1){10.00}}
\put(105.00,25.00){\vector(1,2){5.00}}
\put(110.00,35.00){\vector(-1,1){10.00}}
\put(70.00,35.00){\vector(3,1){30.00}}
\put(70.00,35.00){\vector(-3,1){30.00}}
\put(25.00,35.00){\vector(3,2){15.00}}
\put(5.00,35.00){\vector(1,2){5.00}}
\put(10.00,45.00){\vector(-1,2){5.00}}
\put(10.00,45.00){\vector(3,2){15.00}}
\put(40.00,45.00){\vector(-3,2){15.00}}
\put(25.00,55.00){\vector(-1,2){5.00}}
\put(20.00,65.00){\vector(1,1){10.00}}
\put(5.00,55.00){\vector(3,2){15.00}}
\put(40.00,45.00){\vector(3,1){30.00}}
\put(70.00,55.00){\vector(-3,2){15.00}}
\put(100.00,45.00){\vector(-3,1){30.00}}
\put(100.00,45.00){\vector(1,1){10.00}}
\put(110.00,55.00){\vector(-1,2){5.00}}
\put(120.00,45.00){\vector(-1,1){10.00}}
\put(110.00,55.00){\vector(-3,1){30.00}}
\put(80.00,65.00){\vector(-1,2){5.00}}
\put(70.00,55.00){\vector(1,1){10.00}}
\put(105.00,65.00){\vector(-3,1){30.00}}
\put(25.00,55.00){\vector(3,1){30.00}}
\put(55.00,65.00){\vector(-1,2){5.00}}
\put(55.00,65.00){\vector(1,1){10.00}}
\put(80.00,65.00){\vector(-3,2){15.00}}
\put(65.00,75.00){\vector(-1,2){5.00}}
\put(20.00,65.00){\vector(3,1){30.00}}
\put(50.00,75.00){\vector(1,1){10.00}}
\put(75.00,75.00){\vector(-3,2){15.00}}
\end{picture}

\begin{center}           Fig. 4 \end{center}

Fishburn conjectured that {\em the size of any peak-pit CD does
not exceed $\Phi (n)$}, and verified this conjecture for $n\le 6$.\medskip

Galambos and  Reiner \cite{GR} considered a class of CDs, which we call
\emph{GR-domains} (see the definition in the Appendix), and raised a weakened
version of Fishburn's conjecture saying that \emph{the size of any GR-domain
does not exceed $\Phi (n)$}. It should be noted that an equivalent conjecture
in terms of pseudo-line arrangements was raised earlier by Knuth~\cite{Knuth}.

Monjardet~\cite{Monj} calls a CD {\em connected} if it induces a connected
subgraph of the Bruhat graph. He conjectured that \emph{the size of any
connected CD does not exceed $\Phi (n)$}.

Applying Theorem~4, one can conclude that the conjectures by Fishburn, by
Galambos and Reiner, and by Monjardet are equivalent, and we can express this
conjecture as follows:
   \begin{itemize}
\item[(C)] \emph{the maximum possible size $\gamma_n$ of a tiling CD for $n$ is equal
to $\Phi(n)$.}
  \end{itemize}

However, (C) is not true in general. The authors learnt via B.~Monjardet
(however, without pointing out to us any details or references) that Ondjey
Bilka had established some lower bound on $\gamma_n$ which leads to a
contradiction with~(C). Subsequently the authors found a simple argument, as
follows.\medskip

Let $T$ and $T'$ be  rhombus tilings of zonogons $Z_n$ and $Z_{n'}$,
respectively. We identify the set $[n']$ with the subset $\{n+1,\ldots,n+n'\}$
in $[n+n']$ and merge the top vertex $t(T)$ of $T$ with the bottom vertex
$b(T')$ of $T'$ (erecting $T'$ over $T$). This gives a ``partial tiling'' of
the zonogon $Z_{n+n'}$, as illustrated  in Fig.~5 where $n=4$ and $n'=3$.

\unitlength=.6mm \special{em:linewidth 0.4pt}
\linethickness{0.4pt}
\begin{picture}(90.00,85.00)(-60,4)
\put(60.00,5.00){\vector(-3,1){15.00}}
\put(45.00,10.00){\vector(-1,1){10.00}}
\put(35.00,20.00){\vector(-1,3){5.00}}
\put(30.00,35.00){\vector(0,1){15.00}}
\put(30.00,50.00){\vector(1,3){5.00}}
\put(35.00,65.00){\vector(1,1){10.00}}
\put(45.00,75.00){\vector(3,1){15.00}}
\put(60.00,5.00){\vector(3,1){15.00}}
\put(75.00,10.00){\vector(1,1){10.00}}
\put(85.00,20.00){\vector(1,3){5.00}}
\put(90.00,35.00){\vector(0,1){15.00}}
\put(90.00,50.00){\vector(-1,3){5.00}}
\put(85.00,65.00){\vector(-1,1){10.00}}
\put(75.00,75.00){\vector(-3,1){15.00}}
\put(60.00,5.00){\vector(0,1){15.00}}
\put(60.00,20.00){\vector(-1,3){5.00}}
\put(55.00,35.00){\vector(-1,1){10.00}}
\put(45.00,45.00){\vector(-3,1){15.00}}
\put(30.00,50.00){\vector(3,1){15.00}}
\put(45.00,55.00){\vector(1,1){10.00}}
\put(55.00,65.00){\vector(1,3){5.00}}
\put(75.00,10.00){\vector(0,1){15.00}}
\put(75.00,25.00){\vector(-1,3){5.00}}
\put(70.00,40.00){\vector(-1,1){10.00}}
\put(60.00,50.00){\vector(-3,1){15.00}}
\put(85.00,20.00){\vector(0,1){15.00}}
\put(85.00,35.00){\vector(-1,3){5.00}}
\put(80.00,50.00){\vector(-1,1){10.00}}
\put(70.00,60.00){\vector(-3,1){15.00}}
\put(60.00,20.00){\vector(3,1){15.00}}
\put(75.00,25.00){\vector(1,1){10.00}}
\put(85.00,35.00){\vector(1,3){5.00}}
\put(55.00,35.00){\vector(3,1){15.00}}
\put(70.00,40.00){\vector(1,1){10.00}}
\put(80.00,50.00){\vector(1,3){5.00}}
\put(45.00,45.00){\vector(3,1){15.00}}
\put(60.00,50.00){\vector(1,1){10.00}}
\put(70.00,60.00){\vector(1,3){5.00}}
\put(45.00,28.00){\makebox(0,0)[cc]{$T$}}
\put(45.00,65.00){\makebox(0,0)[cc]{$T'$}}
\end{picture}
\begin{center}           Fig. 5 \end{center}

This partial tiling can be extended (in a unique way, in fact) to a complete
rhombus tiling $\widehat{T}$ of the whole  zonogon $Z_{n+n'}$. If $\sigma$ is a
snake of $T$ and $\sigma'$ is a snake of $T'$, then the concatenated path
$\sigma\sigma'$ is a snake of $\widehat{T}$. Thus, we obtain the injective
mapping
  $$
\Sigma(T)\times \Sigma(T') \to \Sigma(\widehat{T}),
  $$
which gives the inequality $\gamma_n\gamma_{n'}\le \gamma_{n+n'}$.

Now take both $T$ and $T'$ to be Fishburn's tilings for $n=n'=21$. Using a
precise formula for $\Phi(n)$ from~\cite{GR}, one can compute that
$\Phi(21)=4.443.896$ and $\Phi(42)=19.156.227.207.750$. Then
$\Phi(21)^2=19.748.211.658.816>\Phi(42)$. Thus, $\Phi(42)<\gamma_{42}$,
contradicting~(C).
\medskip

\textbf{Remark.} ~The above construction can be given in terms of
``concatenating'' corresponding peak-pit domains rather than tilings. So
Fishburn's conjecture can be disproved without appealing to Theorem~4.

      \section{Some reformulations}

Any linear order can be realized as a snake of some rhombus tiling. However,
this need not hold for a pair of linear orders. For example, the linear orders
213 and 312 (which together with 123 and 321 form the CD ~$\mathcal D_3(\gets)$
from Section~\ref{sec:CD}) cannot appear in the same tiling.

Let us say that two linear orders $\sigma $ and $\tau $ are {\em strongly
consistent} if there exists a tiling $T$ such that $\sigma ,\tau \in \Sigma
(T)$. For example, $\sigma $ and $\tau $ are strongly consistent if $\sigma \ll
\tau $ (where $\ll$ is defined in Section~\ref{sec:Bruhat}). Using observations
and results from previous sections, we can demonstrate some useful equivalence
relations.
\medskip

\noindent {\bf Proposition 3.} ~\emph{ Let $\sigma $ and $\tau $ be linear
orders on $[n]$. The following properties are equivalent:}

{\rm(i)} ~\emph{$\sigma $ and $\tau $ are strongly consistent;}

{\rm(ii)} ~\emph{the set-system $Id(\sigma)\cup Id(\tau )$ is  separated;}

{\rm(iii)} ~\emph{for each triple in $[n]$, the restrictions of $\sigma $ and
$\tau $ to this triple simultaneously satisfy either peak conditions or pit
conditions (or both);}

{\rm(iv)} ~$Id(\sigma)\cup Id(\tau)=Id(\sigma\vee \tau)\cup Id(\sigma\wedge
\tau)$ (where $\vee,\wedge$ concern the Bruhat lattice);

{\rm(iv$'$)} ~$Id(\sigma)\cup Id(\tau)\subseteq Id(\sigma\vee \tau)\cup
Id(\sigma\wedge \tau)$. \medskip

\emph{Proof.} ~Properties (i) and (ii) are equivalent by Theorem~3.

Properties (i) and (iii) are equivalent by Theorem~2.

To see that (i) implies (iv), observe that if $\sigma $ and $\tau $ occur in a
tiling $T$, then $\mathcal{S}(\sigma \vee \tau )$ and $\mathcal{S}(\sigma
\wedge \tau )$ are the left and right envelopes of the snakes for $\sigma $ and
$\tau $, respectively. Therefore, any vertex of the snake $\mathcal{S}(\sigma
\vee \tau )$ is a vertex of $\mathcal{S}(\sigma)$ or $\mathcal{S}(\tau )$, and
similarly for $\mathcal{S}(\sigma \wedge \tau )$. Conversely, each vertex of
$\mathcal{S}(\sigma)\cup \Sscr(\tau)$ is a vertex of $\mathcal{S}(\sigma \vee
\tau )$ or $\mathcal{S}(\sigma \wedge \tau )$.

Obviously, (iv) implies (iv$'$). Let us prove the converse. Since $\sigma\wedge
\tau \ll \sigma\vee \tau$, the linear orders $\sigma\wedge \tau$ and
$\sigma\vee \tau$ are strongly consistent. By the equivalence of (i) and (ii),
$Id(\sigma\vee \tau)\cup Id(\sigma\wedge \tau)$ is a separated system. Since
$Id(\sigma)\cup Id(\tau)\subseteq Id(\sigma\vee \tau)\cup Id(\sigma\wedge
\tau)$, the set-system $Id(\sigma)\cup Id(\tau)$ is separated as well. Thus,
we obtain (ii), whence $(iv')\Rightarrow(iv)$.
\hfill$\Box$\medskip

      \section*{Appendix               }

Here we briefly outline approaches of Abello \cite{Ab}, Galambos
and      Reiner \cite{GR}, and Chameni-Nembua \cite{Ch-Nam}, and
an  interrelation between them and our       approach.

\subsection*{      Abello}

Let $\mathcal D$ be a CD. Then there exists a casting $c$ such that $\mathcal D
\subseteq \mathcal D (c)$ (see Proposition 1). Abello applied this fact to a
maximal chain $\mathcal C $ in the Bruhat lattice (it had been known that any
chain is a CD). In this case the casting $c$ is unique (and is a peak-pit
casting), so the domain $\mathcal D(c)$, denoted by $\widehat{{\mathcal C}}$,
is a CD as well. We call such a CD  an {\em A-domain} (abbreviating
\emph{Abello's domain}). Abello shows that an A-domain is a complete CD.

Note that different chains can give the same A-domain. Maximal chains $\mathcal
C $ and $\mathcal C '$ are called {\em equivalent} if the A-domains
$\widehat{{\mathcal C}}$ and $\widehat{{\mathcal C '}}$ coincide. In the end
of~\cite{Ab} Abello gives another characterization of this equivalence. A
maximal chain in the Bruhat lattice can be thought of as a reduced
decomposition (a product of standard transpositions $s_i$, ~$i\in[n-1]$) of the
longest permutation $\omega $. Then two chains are equivalent if one reduced
decomposition can be obtained from the other by a sequence of transformations,
each replacing a decomposition fragment of the form $s_i s_j$ with $|i-j|>1$ by
$s_js_i$. This characterization became a starting point in Galambos and
Reiner's approach.

                         \subsection*{    Galambos and Reiner}

Let {\bf C} be an equivalence class of maximal chains in the Bruhat lattice.
Define $\mathcal D ({\bf C}):=\cup _{\mathcal C \in {\bf C}}\,\mathcal C $
(Galambos and Reiner referred to this domain as consisting of ``permutations
visited by an equivalence class of maximal reduced decompositions''). We call
$\mathcal D ({\bf C})$ a {\em GR-domain}. It is easy to see (and Galambos and
Reiner explicitly mention this) that the GR-domains are exactly the A-domains.
Moreover, they give a direct proof (in Theorems~1 and~2 of \cite{GR}) that a
GR-domain is a complete CD.

To give a more enlightening representation for the equivalence classes of
maximal reduced decompositions, Galambos and Reiner used \emph{arrangements of
pseudo-lines} (cf.~\cite{BFZ}). Permutations (or linear orders) from the domain
$\mathcal D ({\bf C})$ are realized in these arrangements as certain cutpaths
(viz. directed cuts). Although they do not prove explicitly that the set of
cutpaths of an arrangement forms a complete CD, this can be done rather easily.
Using a relationship between pseudo-line arrangements and rhombus tilings
(cf.~\cite{El}), one can conclude that the GR-domains (as well as the
A-domains) are exactly CDs of tiling type.

\subsection*{ Chameni-Nembua}

One more interesting approach was proposed by Chameni-Nembua. A sublattice
$\mathcal D$ in the Bruhat lattice is called {\em covering} if the cover
relation in this sublattice is induced by the cover relation in the Bruhat
lattice.

Chameni-Nembua shows that a distributive covering sublattice in the Bruhat
lattice is a CD. Suppose that $\mathcal{D}$ is a maximal distributive covering
sublattice. One can easily see that it contains $\alpha$ and $\omega$, and
hence it contains a maximal chain. Therefore, it is a subset of a unique tiling
CD. On the other hand, since any tiling CD forms a distributive covering
sublattice (see Section~\ref{sec:tiling-CD}), one can conclude that $\mathcal
D$ coincides with this tiling CD.

Thus, Chameni-Nembua's approach gives the same class of CDs as the one of
rhombus tilings.\medskip

\textbf{Acknowledgements.} We thank the anonymous referees for remarks and
useful suggestions intended to improve the presentation stylistically. This
research was supported by RFBR grant 10-01-9311-CNRSL\_\,a.


\begin{thebibliography}{99}

\bibitem{Ab} J.M.~Abello, The weak Bruhat order on $S_n$,
consistent sets, and Catalan numbers, \textsl{SIAM J. on Discrete Math.}
\textbf{4} (1991) 1--16.

\bibitem{BFZ} A.~Berenstein, S.~Fomin, and  A.~Zelevinsky,
Parametrizations of canonical bases and totally positive matrices, {\sl
Adv.~Math.} {\bf 122} (1996) 49–-149.

\bibitem{dB} N.G.~de~Bruijn, Dualization of multigrids, {\sl
J.~Phys.~France} {\bf 47} (1986) 3--9.

\bibitem{Ch-Nam} C.~Chameni-Nembua, R\`egle majoritaire et
distributivit\'e dans le permuto\`edre,  \textsl{Mathematiques Informatique et
Sciences humaines} \textbf{108} (1989) 5--22.

\bibitem{DKK-08} V.I.~Danilov, A.V.~Karzanov and G.A.~Koshevoy, On bases of tropical
Pl\"ucker functions, \textsl{ArXiv}:0712.3996v2[mathCO], 2007.

\bibitem{DKK-10} V.~Danilov, A.~Karzanov and G.~Koshevoy, Pl\"ucker environments,
wiring and tiling diagrams, and weakly separated set-systems, \textsl{Adv.~Math.}
\textbf{224} (2010) 1--44.

\bibitem{El} S.~Elnitsky, Rhombic tilings of polygons and clases of
reduced words in Coxeter groups, {\sl J.~Comb. Theory}, Ser.~A, {\bf 77} (1997)
193--221.

\bibitem{Fish} P.~Fishburn, Acyclic sets of linear orders,
 \textsl{Social Choice and Welfare} \textbf{14} (1997) 113--124.

\bibitem{EG} C.Greene and P.Edelman, Combinatorial correspondences
for Young tableaux, balanced tableaux and maximal chains in the
weak Bruhat order of $S_n$, {\em in} ``Combinatorics and Algebra''
(Proceedings, Boulder Conferences) (C.Greene, Ed.), Contemporary
Mathematics, Vol. 34, AMS, Providence, R.I., 1985

\bibitem{Good97} J.E. Goodman,  Pseudoline arrangements, in:
\textsl{Handbook of Discrete and Computational Geometry}, Goodman and O'Rourke
eds., CRC Press, 1997, pp.~83--110.

\bibitem{GR} A.~Galambos and V.~Reiner, Acyclic sets of linear
orders via the Bruhat Order,  \textsl{Social Choice and Welfare} \textbf{30}
(2008) 245--264.

\bibitem{KS} R.~Kenyon and J.-M.~Schlenker, Rhombic embeddings of planar graphs
with faces of degree 4, {\sl ArXiv}:math-ph/0305057, 2003.

\bibitem{Knuth} D. E. Knuth, Axions and Hulls, \textsl{Lect. Notes
Comput. Sci.}, vol.~\textbf{606}, Springer-Verlag, 1992.

\bibitem{LZ} B.~Leclerc and A.~Zelevinsky, Quasicommuting families of
quantum Pl\"ucker coordinates, \textsl{Amer. Math. Soc. Trans., Ser.~2,}
\textbf{181} (1998) 85--108.


\bibitem{Monj} B.~Monjardet, Acyclic domains of linear orders: a
survey, in: \textsl{The Mathematics of Preference, Choice and Order} (S.~Brams,
W.~Gehrlein, and F.~Roberts, eds.), Springer, 2009, pp. 136--160.

\bibitem{Zig} G.M.~Ziegler, Higher Bruhat orders and cyclic
hyperplane arrangements. {\em Topology} \textbf{32} (1993)
259--279.

\end{thebibliography}
      \end{document}